\documentclass[aps,prl,superscriptaddress,floatfix,nobibnotes]{revtex4}
\usepackage{epsfig}
\usepackage{epstopdf}
\usepackage{color}
\usepackage{float}

\usepackage{graphicx}
\usepackage{longtable}
\usepackage{CJK}
\usepackage{color}
\newcommand{\blue}[1]{{ \color{blue} #1 }}

\ifx\pdfoutput\undefined
\usepackage[dvips]{graphicx}
\usepackage[dvipdfm,
  pdfstartview=FitH,
     CJKbookmarks=true,
          bookmarksnumbered=true,
          bookmarksopen=true,
         colorlinks=true,
          pdfborder=002,
          citecolor=blue%
           ]{hyperref}
\AtBeginDvi{} 
\else
\usepackage[
          pdfstartview=FitH,
          CJKbookmarks=true,
          bookmarksnumbered=true,
           bookmarksopen=ture,
           colorlinks=true,   
           citecolor=blue%
           ]{hyperref}

\begin{document}

\title{ A polynomial approach to the Collatz conjecture}

\author{Feng Pan}
\affiliation{Department of Physics, Liaoning Normal University,
Dalian 116029, China}\affiliation{Department of Physics and
Astronomy, Louisiana State University, Baton Rouge, LA 70803-4001,
USA}

\author{Jerry P. Draayer}
\affiliation{Department of Physics and Astronomy, Louisiana State
University, Baton Rouge, LA 70803-4001, USA}
\date{\today}

\begin{abstract}

The Collatz conjecture is explored
using polynomials based on a binary numeral system.
It is shown that the degree of the polynomials, on average,
decreases after a finite number of steps of the Collatz operation,
which provides a weak proof of the conjecture by {using
induction} with respect to the degree of the polynomials.

\end{abstract}

\maketitle

The Collatz  $3x+1$ problem \cite{[1],[2],[3]} concerns
consecutive Collatz operations $C$
to a given integer $n$ with $C[n]=(3n+1)/2$ if $n$ is odd
and $C[n]=n/2$ if $n$ is even.
The conjecture asserts that there is always a finite
number $k$ of the Collatz operations, after which
$C[C[\cdots C[n]\cdots ]=1$.
In this letter, we only consider odd integers $n$ with
the operation $C_{q}[n]=(3n+1)/2^q$, where $q$
is the largest positive integer with $C_{q}[n]=1$ (mod 2).
The conjecture has been verified to be true
for all $n < 20\times  2^{58}$~\cite{[4]} .
\vskip .3cm
It is well known that a natural number less than $2^{p+1}$
can be expressed in {a}
binary numeral system with
\begin{equation}
n=\sum_{\mu=0}^{p}c_{\mu}2^{\mu},
\end{equation}
where there is a unique positive integer $p$ and a set of
coefficients $\{c_{p},c_{p-1},\cdots,c_{0}\}$ with
$c_{i}=0$ or $1$ for the given $n$.
The sequence of bits, $\{c_{p},c_{p-1},\cdots,c_{0}\}$, is just the
binary representation of the integer $n$.
In the following, we always assume that $p$ is a positive finite integer.
In order to visually realize the Collatz operation on odd integers,
we introduce the polynomial
\begin{equation}\label{2}
F_{p}(x)=x^p+c_{p-1}x^{p-1}+c_{p-2}x^{p-2}+\cdots+c_{2}x^{2}+c_{1}x+1
\end{equation}
of degree $p\geq1$, where $x\equiv2$ is assumed throughout, which corresponds to an odd integer $2^{p}+1\leq n\leq
2^{p+1}-1$ given in (1)
with $n\equiv F_{p}(x)$. Arithmetic operations of the polynomials, such as addition and subtraction
$F_{p_{1}}(x)\pm F_{p_{2}}(x)$,  multiplication $F_{p_{1}}(x)F_{p_{2}}(x)$,
and division  $F_{p_{1}}(x)/F_{p_{2}}(x)$
are defined as usual, for which one only needs to keep in mind that $x\equiv2$,
so
\begin{equation}x-1\equiv 1,~x+1\equiv x^2-1,~2\,x^{t-1}\equiv x^{t}\end{equation}
with $t\geq1$ according to the rules of arithmetic operations on
integers in the binary numeral system.
Hence, the resultants of $F_{p_{1}}(x)\pm F_{p_{2}}(x)$ or $F_{p_{1}}(x)F_{p_{2}}(x)$
are still a polynomial of the same type. Table \ref{tb1}
provides $(x+1)^{q}$ for $q\leq 10$ as examples computed in this way.
The degree $p$ of the polynomial $(x+1)^{q}$ increases with $q$ linearly
as $p=u(q)\equiv {\rm Int}[-1/2 + q\,\ln[3]/\ln[2]]$ for $q\geq1$, where
${\rm Int}[r]$ is the nearest integer of $r$.
\vskip .3cm

\begin{table}[tb]
\caption{The polynomial $(1+x)^{q}$ for $q\leq 10$.}
\begin{center}\footnotesize
\begin{tabular*}{120mm}{c@{\extracolsep{\fill}}ccccc}
\hline
$q$~&The degree of $(1+x)^{q}$  & $(1+x)^{q}$\\
\hline
0~ &0 &1\\
1~ &1   &$x+1$\\
2~ &3  &$x^3+1$\\
3~ &4 &$x^4+x^3+x+1$\\
4~ &6   &$x^6+x^4+1$\\
5~  &7 &$x^7+ x^6+ x^5 + x^4+ x +1 $\\
6~ &9 &$x^9+ x^7+ x^6 +  x^4  +x^3+  1  $\\
7~ &11 &$x^{11}+ x^7 + x^3  +x   +  1$\\
8~ &12 &$x^{12}  + x^{11} + x^8 + x^7 + x^5  +1 $\\
9~ &14 &$ x^{14}+ x^{11} +x^{10} +x^7 +  x^6 +x^5 +   x + 1  $\\
10~&15 &$x^{15}+  x^{14} +x^{13} +  x^{10} +x^9 +  x^7 +x^5 +  x^3 +1$\\
\hline\hline
\end{tabular*}
\end{center}\label{tb1}
\end{table}

The Collatz operation on the polynomial $F_{p}(x)$ representing
an odd integer $2^{p}+1\leq n\equiv F_{p}(x)\leq 2^{p+1}-1$  is defined as
\begin{equation}\label{4}
C_{q}[F_{p}(x)]=x^{-q}\left((x+1)F_{p}(x)+1\right),
\end{equation}
where the positive integer $q\geq1$ is chosen to be the largest
such that $C_{q}[F_{p}(x)]$ is still a polynomial
of the same type defined by (\ref{2}), or equivalently  $(x+1)F_{p}(x)+1$
is factorizable as $(x+1)F_{p}(x)+1=x^{q}\,Q_{t}(x)$,
where $t$ is an positive integer, and
$C_{q}[F_{p}(x)]=Q_{t}(x)$ is a polynomial of the same type
of degree $t$ representing another
odd integer.
Thus, the polynomial ring $R=\{\{F_{0}(x)\},\cdots,\{F_{p}(x)\},\cdots\}$ constructed by
a series of polynomials $\{F_{0}(x)\}$,$\cdots$,$\{F_{p}(x)\}$, $\cdots$,
where $F_{p}(x)$ is given by (\ref{2}), and $\{F_{p}(x)\}$ are formed by
$2^{p-1}$ different combinations of $0$ and $1$
in the sequence of bits $\{c_{p-1},\cdots,c_{1}\}$ in $F_{p}(x)$,
is algebraically closed under the Collatz operation.
The Collatz conjecture can then be stated as follows:
\vskip .3cm
\noindent {\bf The Collatz Conjecture}
Any degree $p$ polynomial $F_{p}(x)\in R$
decreases after a finite number of the Collatz operations.
Namely, $C_{q_{l}}[\cdots[C_{q_{1}}[F_{p}(x)]\cdots]=F_{p^{\prime}}(x)\in R~\forall~F_{p}(x)\in R$
with $p^{\prime}\leq p-1$ and finite $l$,
and eventually  $C_{q_{k}}[\cdots[C_{q_{1}}[F_{p}(x)]\cdots]=F_{0}(x)=1$ with a finite $k\geq l$.

\vskip .3cm
\noindent {\bf Definition} If a polynomial $F_{p}(x)$
defined by (\ref{2}) satisfies $C_{q_{l}}[C_{q_{l-1}}[\cdots C_{q_{1}}[F_{p}(x)]\cdots ]=1$ with finite $l$,
$F_{p}(x)$ is called the Collatz polynomial.
\vskip .3cm {Concerning} the Collatz operation on $R$,
we have the following corollaries:
\vskip .3cm
\noindent {\bf Corollary 1}
If $Q_{t}(x)=C_{q}[F_{p}(x)]$,
then $Q_{t}(x)=C_{q+2}[1+x^2 F_{p}(x)]=C_{q+4}[1+x^2+x^4 F_{p}(x)]=\cdots$, which can be verified directly by
the Collatz operation.
\vskip .3cm
\noindent{\bf Corollary 2} {After} consecutive iterations according to
Corollary 1, the following series $\{U_{k}(x)=\sum_{t=0}^{k}x^{2t}\}$ ($k=0,1,2,\cdots$) satisfies $C_{2k+2}[U_{k}(x)]=U_{0}=1$,
which was {noted} in \cite{[5]}.

\vskip .3cm

If there are $m+2$ nonzero terms in $F_{p}(x)$,
without loss of generality, $F_{p}(x)$ can be expressed as
\begin{equation}\label{5}
F^{(m)}_{p}(x)=x^{p}+\sum_{i=1}^{m}x^{k_{i}}+1,
\end{equation}
in which $p-1\geq k_{m}> k_{m-1}>\cdots> k_{1}\geq 1$
for $1\leq m\leq p-1$ is assumed, or
$F^{(m)}_{p}(x)=x^{p}+1$ with $m=0$.
\vskip .3cm
For a given polynomial $F^{(m)}_{p}(x)$,
the degree of $C_{q_{l}q_{l-1}\cdots q_{1}}[F^{(m)}_{p}(x)]$,
where  $C_{q_{l}q_{l-1}\cdots q_{1}}[F^{(m)}_{p}(x)]\equiv C_{q_{l}}[\cdots
C_{q_{1}}[F^{(m)}_{p}(x)]$ stands for $l$ times
of the Collatz operation on $F^{(m)}_{p}(x)$,
can be expressed as
\begin{equation}\label{42}
{\rm Deg}(C_{q_{l}q_{l-1}\cdots q_{1}}[F^{(m)}_{p}(x)])= p+u(l)+1-\sum_{i=1}^{l}q_{i},
\end{equation}
where
$u(l)= {\rm Int}[-1/2 + l\,\ln[3]/\ln[2]]$ is the degree of $(x+1)^l$,
which is  obtained simply based on the power counting of the leading term of
$C_{q_{l}q_{l-1}\cdots q_{1}}[F^{(m)}_{p}(x)]$.
Concerning the Collatz operation, we have the following proposition:
\vskip .3cm
\noindent{\bf Proposition 1} Except {for the} trivial polynomial $F_{0}=U_{0}=1$,
{the Collatz operation on $F_{p}(x)$, $C_{q}[F_{p}(x)]$, is different from $F_{p}(x)$ itself.}
\vskip .3cm
The validity of Proposition 1 is obvious from the Collatz operation
on $F_{p}(x)$. If there is  an $F_{p}(x)$ satisfies $C_{q}[F_{p}(x)]=F_{p}(x)$,
one can deduce that $F_{p}(x)=1/(x^{q}-x-1)$. Since $q\geq 1$ and $F_{p}(x)>1$ being
a polynomial,
the only possible solution for  $F_{p}(x)\geq1$
is $F_{p}(x)=1$ with $q=2$,
which is excluded in Proposition 1.
Directly taking the Collatz operation on $F_{p}(x)$,
we also have
\vskip .3cm
\noindent{\bf Corollary 3} The degree of the
polynomials $F^{(m)}_{p}(x)$ for $m=0$, $m=1$,
and $m\geq 1$ with $k_{1}\geq 2$ and $k_{m}<p-1$
decreases after a few steps of  the Collatz operation.
\vskip .3cm
For $m=0$, we have
\begin{eqnarray}\label{61}
&C_{2}[F^{(0)}_{p}(x)]]=x^{p-1}+x^{p-2}+1~{\rm for}~ p\geq2,~~
C_{4}[C_{1}[F^{(0)}_{p}(x)]]=1~{\rm for}~p=1.
\end{eqnarray}
For $m=1$ and $p-1\leq k_{1}\leq 1$, the resultants of direct Collatz operations on
$F^{(1)}_{p}(x)$  are
\begin{eqnarray}\label{6}\nonumber
&C_{4}[C_{1}[F^{(1)}_{p}(x)]]=x^{p-2}+x^{p-5}+1~{\rm for}~ k_{1}=1,\\\nonumber
&C_{2}[F^{(1)}_{p}(x)]=x^{p-1}+x^{p-2}+x^{k-1}+x^{k-2}+1~{\rm for}~p-2\geq k_{1}\geq2,\\
&C_{2}[C_{2}[F^{(1)}_{p}(x)]]=x^{p-1}+x^{p-2}+x^{p-4}+x^{p-5}+1~{\rm for}~k_{1}=p-1.
\end{eqnarray}
For $m\geq 1$ and $k_{1}\geq2$ and $k_{m}<p-1$,
\begin{eqnarray}\label{610}
C_{2}[F^{(m)}_{p}(x)]=x^{p-1}+x^{p-2}+\sum_{i=1}^{m}x^{k_{i}-1}+\sum_{i=1}^{m}x^{k_{i}-2}+1.
\end{eqnarray}
It is obvious that the degrees of the resultants shown in (\ref{61}), (\ref{6}), and (\ref{610})
are less than $p$.

\vskip .3cm
Furthermore, for $p-1\leq k_{2}<k_{1}\leq1$, the resultants of direct Collatz operations on
$F^{(2)}_{p}(x)$  are
\begin{eqnarray}\label{7}\nonumber
&C_{1}[C_{1}[F^{(2)}_{p}(x)]]=x^{p+1}+x^{p-2}+x^4+1~{\rm for}~ k_{1}=1,~k_{2}=2,\\\nonumber
&C_{1}[F^{(2)}_{p}(x)]=x^{p}+x^{p-1}+x^{4}+1~{\rm for}~k_{1}=1,~k_{2}=3,\\\nonumber
&C_{3}[C_{1}[F^{(2)}_{p}(x)]]=x^{p-1}+x^{p-4}+x^{3}+x+1~{\rm for}~k_{1}=1,~k_{2}=4,\\\nonumber
&C_{4}[C_{1}[F^{(2)}_{p}(x)]]=x^{p-2}+x^{p-5}+x^{k_{2}-2}+x^{k_{2}-5}+1~{\rm for}~k_{1}=1,~k_{2}>4,\\\nonumber
&C_{3}[F^{(2)}_{p}(x)]=x^{p-2}+x^{p-3}+x^2+1~{\rm for}~ k_{1}=2,~k_{2}=3,\\\nonumber
&C_{4}[F^{(2)}_{p}(x)]=x^{p-3}+x^{p-4}+x^{k_{2}-3}+x^{k_{2}-4}+1~{\rm for}~ k_{1}=2,~k_{2}\geq4,\\
&C_{2}[F^{(2)}_{p}(x)]=x^{p-1}+x^{p-2}+x^{k_{2}-1}+x^{k_{2}-2}+x^{k_{1}-1}+x^{k_{1}-2}+
1~{\rm for}~ k_{1}\geq3.
\end{eqnarray}
The above examples show that the first two cases of $F^{(2)}_{p}(x)$ given in (\ref{7})
increase or {remain unchanged in their degree}
after a few steps of the Collatz operation, while {the other}
cases decrease in the degree,
which are mainly determined
by the values of $k_{i}$ ($i=1,\cdots,m$),
especially by those of  $k_{m}$ and $k_{1}$ of the monomials $x^{k_{m}}$ and $x^{k_{1}}$ involved.
When $p-k_{m}=1$, the degree of the leading term of $(x+1)(x^{p}+x^{k_{m}})$ will
increase from $p$ to ${p+2}$, while $(x+1)(x^{k_{1}}+1)+1$ becomes $x^{3}+x$
when $k_{1}=1$, with which the divisor $x^{q}$ required in the Collatz operation is
the smallest in $q$ with $q=1$
and resulting in  $p+1$  degree polynomial $C_{1}[F^{(m)}_{p}(x)]$ for $m>1$.
This situation {remains} unchanged after several steps of the Collatz
operation, especially when $k_{m}=p-1$, $k_{m-1}=p-2$, $\cdots$, $k_{1}=1$ for $m=p-1$.
As the consequence, the worst situation is typically represented by $F^{(m)}_{p}(x)=x^{p+1}-1$
with
$m=p-1$, which needs most steps of the Collatz operation
to get a polynomial with degree less than $p$.
Especially, $\sum_{i=1}^{l}q_{i}=l$ for $l\leq p$ for this case, with
which the degree of the resulting polynomial will increase within
the first $p$ steps of the Collatz operation.
Table 2\ref{tb2} provides the resultant of
$F^{(p-1)}_{p}(x)$ after $p$ steps of the Collatz operation
for $p\leq 32$ explicitly as examples, in which the last column
provides the $\sum_{i=1}^{k}q_{i}/k$, where $k$ is the total number of steps of the Collatz operation
needed for $C_{1}^{(p)}[F^{(p-1)}_{p}(x)]$ to reach $F_{0}$.
In this case, one can verify that
\begin{equation}\label{8}
G_{u(p)+1}(x)=C_{1}^{(p)}[F^{(p-1)}_{p}(x)]=
({x+1\over{x}})^{p}F^{(p-1)}_{p}(x)+{1\over{x}}
\sum_{\mu=0}^{p}({x+1\over{x}})^{\mu}=
x(x+1)^{p}-1,
\end{equation}
where $C_{1}^{(p)}$ stands for $p$ times of the Collatz operation $C_{1}$,
which is a polynomial of degree $u(p)+1=1+{\rm Int}[-1/2+p \ln[3]/\ln[2]]>p$.
Obviously,
$q_{i}=1$ for $i\leq p$ within the first $p$ steps of the Collatz operation
on $F_{p}^{(p-1)}(x)$.

\newpage

{\noindent {\bf Table 2}. Some $F^{(p-1)}_{p}(x)$ after
$p$ steps of Collatz operations $G_{u(p)+1}(x)=C_{1}^{(p)}[F^{(p-1)}_{p}(x)]$,
where, due to (\ref{12}), only  $G_{u(p)+1}(x)$ with even $p$ for $p\leq 32$  are provided.}
\begin{center}\footnotesize
\begin{tabular*}{139mm}{c@{\extracolsep{\fill}}ccccc}
\hline
$p$~ & $u(p)+1$ & $G_{u(p)+1}(x)$ &$\langle \sum_{i=1}^{k}q_{i}\rangle/k$\\
\hline
2~ &4 &$1+x^4$ &3\\
4~ &7 &$1+x^5+x^7$&2.172\\
6~ &10&$1+x^4+x^5+x^7+x^8+x^{10}$ &2.778\\
8~ &13 &$1 + x^6 + x^8 + x^9 + x^{12} + x^{13}$&2.357\\
10~ &16 &$1 + x^4 + x^6 + x^8 + x^{10} + x^{11} + x^{14} + x^{15} + x^{16}$&1.957\\
12~ &20 &$1 + x^5 + x^6 + x^7 + x^8 + x^9 + x^{10} + x^{12} + x^{13} + x^{20}$&2.045\\
14~ &23 &$1 + x^4 + x^5 + x^6 + x^7 + x^9 + x^{10} + x^{12} + x^{13} + x^{14} + x^{15} +x^{16} + x^{20} + x^{23}$&3.033\\
16~ &26 &$1 + x^7 + x^9 + x^{10} + x^{11} + x^{13} + x^{15} + x^{16} + x^{21} + x^{24} + x^{26}$&1.968\\
18~ &29 &$1 + x^4 + x^7 + x^9 + x^{13} + x^{16} + x^{17} + x^{18} + x^{19} + x^{21} + x^{25}
+ x^{26} + x^{27} + x^{29}$&2.333 \\
20~ &32 &$1 + x^5 + x^8 + x^9 + x^{10} + x^{12} + x^{13} + x^{19} + x^{21} + x^{23} + x^{24}+$\\
&&$ x^{25} + x^{26} + x^{27} + x^{28} + x^{31} + x^{32}$&2.072\\
22~ &35 &$1 + x^4 + x^5 + x^{12} + x^{13} + x^{14} + x^{15} + x^{16} + x^{19} + x^{21} + x^{22} + x^{23}+$\\
&&$ x^{26} + x^{27} + x^{28} + x^{31} + x^{33} + x^{34} + x^{35}$&1.822\\
24~ &39 &$1 + x^6 + x^7 + x^8 + x^{12} + x^{13} + x^{14} + x^{16} + x^{19} + x^{20} + x^{21}+$\\
&&$x^{26} + x^{31} + x^{32} + x^{33} + x^{39}$ &1.985\\
26~ &42 &$1 + x^4 + x^6 + x^7 + x^8 + x^9 + x^{10} + x^{11} + x^{12} + x^{13} + x^{14} +x^{15} + x^{18} +$\\
&&$x^{25} + x^{26} + x^{29} + x^{31} + x^{32} + x^{33} + x^{34} + x^{35} + x^{36} + x^{39} + x^{42}$&1.955\\
28~ &45 &$1 + x^5 + x^6 + x^9 + x^{10} + x^{11} + x^{12} + x^{13} + x^{14} + x^{15} + x^{18}+ x^{19} +x^{21} + $\\
&&$x^{25} + x^{26} + x^{28} + x^{30} + x^{31} + x^{34} + x^{35} + x^{36}+ x^{39} + x^{40} + x^{43} + x^{45}$
&2.040\\
30~&48 &$1 + x^4 + x^5 + x^6 + x^8 + x^{12} + x^{13} + x^{14} + x^{15} + x^{18} + x^{20}+x^{23} + x^{24} +$\\
&&$x^{25} + x^{26} + x^{31} + x^{32} + x^{33} + x^{39} + x^{41} + x^{42} +x^{44} + x^{45} + x^{46} + x^{48}$&2.048\\

32~ &51&$1 + x^8 + x^{10} + x^{11} + x^{12} + x^{13} + x^{14} + x^{18} + x^{19} + x^{20}+x^{21} + x^{26} + x^{30} +$\\
& &$ x^{31} + x^{32} + x^{33} + x^{34} + x^{35} + x^{36} + x^{39}+
x^{41} + x^{43} + x^{45} + x^{48} + x^{50} + x^{51}$ &1.924\\
\hline\hline
\end{tabular*}\label{tb2}
\end{center}

\vskip .3cm
From (\ref{8}), we have
\begin{eqnarray}\label{9}
&G_{u(p)+1}(x)=x(x+1)^{p}-1,~~ G_{u(p+1)+1}(x)=(x+1)G_{u(p)+1}(x)+x,
~~C_{q}[G_{u(p)+1}(x)]={(x+1)^{p+1}-1\over{x^{q-1}}},
\end{eqnarray}
where $q>2$ for $p$ odd, and $q=2$ for $p$ even.
Hence,
\begin{equation}\label{10}
C_{r}[C_{2}[G_{u(p)+1}(x)]]={(x+1)^{p+2}-1\over{x^{r+1}}}
\end{equation}
for $p$ even, where $r\geq 2$ because $p+1$ is odd.
Moreover,

\begin{equation}\label{11}
C_{r+2}[G_{u(p+1)+1}(x)]={(x+1)^{p+2}-1\over{x^{r+1}}}.
\end{equation}
Combining Eqs. (\ref{10}) and (\ref{11}), we get
\begin{equation}\label{12}
C_{r+2}[G_{u(p+1)+1}(x)]=C_{r}[C_{2}[G_{u(p)+1}(x)]]~~{\rm for~}~p=0,\,2,\,4,\,\cdots.
\end{equation}
According to (\ref{12}), if  $G_{u(p)+1}(x)$ is a Collatz polynomial,
$G_{u(p+1)+1}(x)$ is also a Collatz polynomial, which, however, is valid
only when $p$ is even.
In addition, according to (\ref{9}) and (\ref{10}),
\begin{equation}\label{13}
C_{r}[C_{2}[G_{u(p)+1}(x)]]={(x+1)^{p+2}-1\over{x^{r+1}}}=
{x(x+1)^{p+2}-x\over{x^{r+2}}}={x(x+1)^{p+2}-1-1\over{x^{r+2}}}
={G_{u(p+2)+1}(x)-1\over{x^{r+2}}},
\end{equation}
from which we have
\begin{equation}\label{14}
G_{u(p+2)+1}(x)=x^{r+2}C_{r}[C_{2}[G_{u(p)+1}(x)]]+1
\end{equation}
for $p=0,\,2,\,4,\,\cdots$, where the positive integer
$r$ changes with $p$ quasi-periodically,
of which some examples are provided in Table 3 and consistent with
the results shown in Table 2.

\newpage
{\noindent {\bf Table 3}. The positive integer $r$  in Eq. (\ref{14})
as a function of even $p$ for $p\leq 32$.}
\begin{center}\footnotesize
\begin{tabular*}{139mm}{c@{\extracolsep{\fill}}ccccc}
\hline
$p$~  &$r$ &$p$~  &$r$\\
\hline
0 &2 &$x$~  &3\\
$x^2$ &2 &$x^2+x$ &4\\
$x^3$  &2 &$x^3+x$ &3  \\
$x^3+x^2$   &2 &$x^3+x^2+x$ &5\\
$x^4$   &2  &$x^4+x$ &3  \\
$x^4+x^2$   &2 &$x^4+x^2+x$ &4\\
$x^4+x^3$   &2 &$x^4+x^3+x$ &3  \\
$x^4+x^3+x^2$   &2 &$x^4+x^3+x^2+x$ &6\\
$x^5$   &2\\

\hline\hline
\end{tabular*}\label{tb3}
\end{center}

\vskip .3cm
In order to demonstrate the patten of $\{F_{p}(x)\}$
after the Collatz operation, one may denote $F_{p}(x)$ and
the resultant of $F_{p}(x)$ after several steps of the Collatz operation
as nodes. If  $F_{p'}(x)$
is the resultant of  $C_{q}[F_{p}(x)]$,
$F_{p}(x)$ and $F_{p'}(x)$ are connected with
an arrow line, of which the arrow points to
$F_{p'}(x)$. Hence, one can generate the Collatz
tree graph for  $\{F_{p}(x)\}$ under the
Collatz operation.

\begin{figure}[H]
\begin{center}
  \includegraphics[width=5.in]{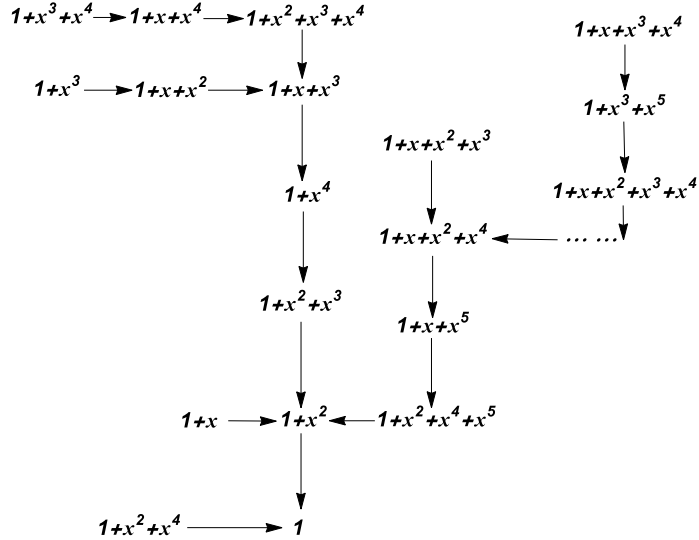}
  \end{center}
  \caption{A part of the Collatz tree graph for the polynomials $\{F_{p}(x)\}$ with $p\leq 4$. \label{fig1}}
\end{figure}

Fig. \ref{fig1} shows a part of the Collatz tree graph for $\{F_{p}(x)\}$
with $p\leq 4$, where a long path with $29$ nodes from
$F_{4}^{(3)}$ to $1+x+x^2+x^4$, in which
no polynomial $F_{m}(x)$ with $m\leq 4$ appears, is
abbreviated with dots.
Though Fig 1 only provides with
a small portion of the graph with $p\leq 4$, its patten
is quite the same as that of the whole tree graph due to the fact
that the properties of the polynomials $\{F_{p}(x)\}$
for either even $p$  or those for  odd $p$ are the same among themselves.
The common features of the graph can be summarized as follows:
{\bf (a)} Due to Proposition 1, the only endpoint node
on the tree graph is $F_{0}$, and there is no other
endpoint node on any path. Moreover, Proposition 1
also asserts that the tree is unique.
Namely, there is no other separate trees containing
some of the polynomials with a different endpoint.
{\bf (b)} For a given $F_{p}(x)$, there is one and only one path on the tree
graph,
which can easily be proven because
$C_{q}[F_{p}(x)]$ is unique.
This unique path is
towards $F_{0}$, which is due to
the fact that $F_{0}$ is the only endpoint node on the tree graph.
Actually, $F_{0}$ is the unique invariant polynomial
under the Collatz operation.
{\bf (c)} If $F_{p}(x)$ is a
factor of $x+1\equiv x^2-1$, it is a starting node of a path
to $F_{0}$, because there is no polynomial after the
Collatz operation to be a factor of $x+1\equiv x^2-1$.
{\bf (d)} From $F_{p}(x)$ towards $F_{0}$, in comparison with
the former node on the move indicated by the arrow,
the degree of the polynomials on the path
may keep unchanged, decrease, or increase.
Without exception, increasing in the degree
with amount $\Delta p \geq2$ occurs only
on the path from $F_{m}^{(m-1)}(x)$ towards
$F_{0}$ for $m\leq p$, which shows that
$F_{m}^{(m-1)}(x)$ with $m\leq p$ are the only sources
resulting in far increasing of the degree of the polynomial
$F_{p}(x)$ after the Collatz operation.

\vskip .3cm
{Concerning} the common features of the tree graph,
one of the unsolved problems is {the possibility of}
circles under the Collatz operation.
Let $\{f_{k}\}$ ($k=1,2,\cdots, n$) are a set of
odd integers satisfying
$f_{\mu+1}=C_{q_{\mu}}[f_{\mu}]$ for $\mu=1,2,\cdots, k-1$.
If $f_{1}=C_{q_{k}}[f_{k}]$, then  $\{f_{k}\}$ forms
a circle \cite{[7]} under the Collatz operation.
Another unsolved problem is that there
may be some polynomials $\{g_{k}\}$, of
which the degrees go to infinity after
consecutive Collatz operations.
In these two cases, the related polynomials will
be separated from the tree graph shown in Fig. \ref{fig1}.
It was shown in \cite{[7]} that the circle containing $F_{p}(x)$
is not possible for small $p$, but might occur when $p$ is large.
For both cases, the degree of  $F_{p}(x)$, {on average,}
will never decrease under the Collatz operation.
However, the structure of the polynomials
are all the same. Namely, whether a polynomial
$F_{p}(x)$ is the Collatz polynomial should be
independent of the degree $p$.
A short discussion on this problem will be
made later on.

\vskip .3cm
Moreover,  when $m=2k$, $F_{2k}^{(2k-1)}(x)$
can be expressed as
\begin{equation}\label{15}
F_{2k}^{(2k-1)}(x)=C_{2}[\sum_{\mu=1}^{k+1}x^{2\mu-1}-1],
\end{equation}
where $\sum_{\mu=1}^{k+1}x^{2\mu-1}-1$ is a polynomial
of degree $2k+1$. While $F_{m}^{(m-1)}(x)$ is the starting
node of a path when $m$ is odd, because
$F_{2k+1}^{(2k)}$ is a factor of $x+1$.
When $2k=3\tau+1$ for $\tau=1,\,3,\,5\cdots$,

\begin{equation}\label{16}
\sum_{\mu=1}^{k+1}x^{2\mu-1}-1=C_{1}[H_{2k}(x)],
\end{equation}
for $k=2,\,3,\,\cdots$,
where
\begin{eqnarray}\label{17}\nonumber
&H_{4}(x)=1+x+x^3+x^4,\\
&H_{4+6(t+1)}(x)=H_{4+6t}\,(x) +(1+x+x^2)x^{6t+8}
\end{eqnarray}
for $t=0,\,1,\,2,\,\cdots$.
Hence, when $p=2k=3\tau+1$ for $\tau=1,\,3,\,5,\,\cdots$,
$H_{2k}(x)$ and $F_{2k-1}^{(2k)}$ are on the same
path, where
\begin{equation}\label{18}
F^{(2k-1)}_{2k}(x)=C_{2}[C_{1}[H_{2k}(x)]]
\end{equation}
for $k=C_{1}[\tau]$ with $\tau=1,\,3,\,5\,\cdots$.

\vskip .3cm
Therefore, if no circle and non-decreasing in the degree occur under the Collatz operation,
the polynomials $\{F_{p}(x)\}$ under
the consecutive steps of the Collatz operation
like a simple board game, for which the rules of the moves
are determined by the Collatz operation.
Once a move starts from any one of the nodes
on the board, there is only one possible path, which is towards the only destination
$F_{0}=1$ under the Collatz operation.
Since $F_{0}=1$ is the only endpoint node on
the tree graph,
the degree of the polynomials $\{F_{p}(x)\}$
after finite steps of the Collatz operation,
{on average,} should decrease, which
can be estimated as follows:

\vskip .3cm
\noindent{\bf Proposition 2} The degree $p$ of the Collatz polynomials
$\{F_{p}(x)\}$, {on average,} decreases after a finite
steps of the Collatz operation.
\vskip .3cm
Generally, the polynomials (\ref{5}) may  be one of the following four
cases with
\begin{eqnarray}\label{4211}
F_{p}(x)\equiv F^{(m)}_{p}(x)=\left\{
\begin{array}{l}
1+f^{1}_{p}(x),\\
1+x+f^{2}_{p}(x),\\
1+x^2+f^{3}_{p}(x),\\
1+x+x^2+f^{4}_{p}(x)
\end{array}\right.
\end{eqnarray}
for $m\geq 3$,
where $f^{1}_{p}(x)=F^{(m)}_{p}(x)-1$ when $k_{1}\geq 3$ in $F^{(m)}_{p}(x)$,
$f^{2}_{p}(x)=F^{(m)}_{p}(x)-1-x$ when $k_{1}=1$ and $k_{2}\geq 3$ in $F^{(m)}_{p}(x)$,
$f^{3}_{p}(x)=F^{(m)}_{p}(x)-x^{2}-1$ when $k_{1}=2$ and $k_{2}\geq 3$ in $F^{(m)}_{p}(x)$,
and $f^{4}_{p}(x)=F^{(m)}_{p}(x)-x^{2}-x-1$ when $k_{1}=1$, $k_{2}=2$, and $k_{3}\geq 3$ in $F^{(m)}_{p}(x)$.
The Collatz operation $C_{q_{1}}$ to the first case is definitely with $q_{1}=2$,
$C_{q_{2}}$ or $C_{q_{4}}$ to the second or the fourth case is with $q_{2}=q_{4}=1$,
and $C_{q_{3}}$ to the third case may be with  $q_{3}=3$,  $4$, $\cdots$, $p+2$.
Anyway, $q_{i}\geq1$ ($i=1,2,3,4$).
Though it is difficult to calculate the distribution of these $q_{3}$ values
for $F^{(m)}_{p}(x)$,
$q_{3}=3$ is the smallest value, which is most possible to occur in this case.
Due to the fact that
\begin{eqnarray}\label{4212}
\begin{array}{l}
C_{2}[1]=1,\\
C_{1}[1+x]=1+x^2,\\
C_{1}[1+x+x^2]=(1+x)+x^3,
\end{array}\end{eqnarray}
the resultant of $F^{(m)}_{p}(x)$ after  the
Collatz operation, $C_{q_{i}}[F^{(m)}_{p}(x)]$
($i=1,2,3,4$), is still of the form shown in (\ref{4211}).
Though there are more possible outcomes for the third case
of (\ref{4211}), the resultant of $C_{q_{3}}[F^{(m)}_{p}(x)]$
is also of the form shown in (\ref{4211}).
The exceptional case is $F^{(m)}_{p}(x)$
with $m=p-1$, for which $q_{i}=1$ for $i\leq p$.
Except $F^{(p-1)}_{p}(x)$ and some polynomials related to it under the Collatz operation,
the four cases listed in (\ref{4211}) occur approximately
randomly after the Collatz operation
on the polynomials $F_{p}(x)$
including  $G_{u(p)+1}=C_{1}^{(p)}[F^{(p-1)}_{p}(x)]$ shown in (\ref{9})
with some examples provided in Table 2{\ref{tb2}}.
Hence, except $F^{(p-1)}_{p}(x)$ and some polynomials
related to it under the Collatz operation, {on average,} the lower bound of the
mean-value $\langle \sum_{i=1}^{l}q_{i}\rangle$
after $l$ steps of the Collatz operation on $F_{p}(x)$
 can be expressed as

\begin{equation}\label{4200}
\langle \sum_{i=1}^{l}q_{i}\rangle\geq\left(
{2\over{4}}+{1\over{4}}+{3\over{4}}+{1\over{4}}\right)l= 1.75\,l.
\end{equation}
It should be stated that (\ref{4200}) only provides the lower bound
of the mean-value for the polynomial $F_{p}(x)$, which is
not related to $F^{(p-1)}_{p}(x)$  under the Collatz operation.
For example, the values of $\langle \sum_{i=1}^{k}q_{i}\rangle/k$ shown in Table 2
are always greater than $1.75$ given in  (\ref{4200}).
There are more extreme cases with
$\sum_{i=1}^{l}q_{i}$ deterministically greater than $\langle \sum_{i=1}^{l}q_{i}\rangle$.
For example, since the Collatz operation on $\{U_{p}(x)\}$ provided
in Corollary 2 always results in $C_{2p+2}[U_{p}(x)]=U_{0}=1$,
especially $C_{2}[U_{0}(x)]=U_{0}=1$, $q_{1}=2p+2$ for the former
and $\sum_{i=1}^{k}q_{i}=2k$ for the latter.  For a class of polynomials $\{\bar{F}_{p}(x)\}$
related to $F^{(p-1)}_{p}(x)$ under the Collatz operation,
the mean-value $\langle \sum_{i=1}^{l}q_{i}\rangle$ should be modified as
\begin{equation}\label{4201}
\langle \sum_{i=1}^{l}q_{i}\rangle
\left\{
\begin{array}{l}
~~\approx l~~~~~~{\rm for}~l\leq p,\\
\\
\geq 1.75\,l~~~{\rm for}~l> p.
\end{array}\right.
\end{equation}
Thus, it can now be shown that Proposition 2 is valid
in concerning (\ref{42}) and the estimation of the lower bound of the
mean-value $\langle \sum_{i=1}^{l}q_{i}\rangle$
after $l$ steps of the Collatz operation  (\ref{4200}) and (\ref{4201}).
Since, after $l$ steps of the Collatz operation,
the degree of $C_{q_{l}q_{l-1}\cdots q_{1}}[F_{p}(x)]$,

\begin{equation}\label{4202}
{\rm Deg}(C_{q_{l}q_{l-1}\cdots q_{1}}[F_{p}(x)])=
p+u(l)+1-\sum_{i=1}^{l}q_{i},
\end{equation}
and, on the average, $\langle \sum_{i=1}^{l}q_{i}\rangle<\sum_{i=1}^{l}q_{i}$,
\begin{equation}\label{4203}
{\rm Deg}(C_{q_{l}q_{l-1}\cdots q_{1}}[F_{p}(x)])<
p+u(l)+1-\langle \sum_{i=1}^{l}q_{i}\rangle=p+{1\over{2}}+\left(
{\ln[3]\over{\ln[2]}}-1.75\right)\,l\simeq p+0.5- 0.165037\,l
\end{equation}
for $F_{p}(x)$ not related to $F^{(p-1)}_{p}(x)$ under the Collatz operation,
while
\begin{equation}\label{4204}
{\rm Deg}(C_{q_{l}q_{l-1}\cdots q_{1}}[F^{(p-1)}_{p}(x)])\left\{
\begin{array}{l}
\approx p+0.5+ 0.58496\,l~~{\rm for}~~l\leq p,\\
\\
<0.5- 0.165037\,l+1.75\,p~~{\rm for}~~l>p,
\end{array}\right.
\end{equation}
which also applies to $\{\bar{F}_{p}(x)\}$.
(\ref{4203}) and (\ref{4204}) provide with
the average upper bound of
${\rm Deg}(C_{q_{l}q_{l-1}\cdots q_{1}}[F_{p}(x)])$.
It can be observed from (\ref{4203}) and (\ref{4204}) that
${\rm Deg}(C_{q_{l}q_{l-1}\cdots q_{1}}[F_{p}(x)])$
with sufficiently large and finite $l$ will be less than $p$, with which
$C_{q_{l}q_{l-1}\cdots q_{1}}[F_{p}(x)]<F_{p}(x)$ is definitely satisfied.
Hence, the Proposition 2 is stronger than and consistent with the results shown in \cite{[1]}.
In addition the number of steps of the Collatz operation needed
for $F_{p}(x)$ to reach $F_{0}$ estimated by (\ref{4203}) or (\ref{4204})
is slightly larger than that estimated by previous probabilistic prediction
\cite{[6]}, because (\ref{4203}) or (\ref{4204}) provides with
the upper bound of ${\rm Deg}(C_{q_{l}q_{l-1}\cdots q_{1}}[F_{p}(x)])$.
Proposition 2 also ensures that there is
no circle containing $F_{p}(x)$ for any $p$, and
there is no polynomial non-decreasing in its degree
after
consecutive Collatz operations. Thus, $F_{0}$
is the only endpoint of the unique tree graph.

\vskip .3cm

Since $F_{0}$ is the obvious Collatz polynomial,
if $\{F_{\mu}(x)\}$ with $\mu=0,1,\cdots,p-1$
have been verified to be the Collatz polynomials,
$\{F_{p}(x)\}$ are also the Collatz polynomials
because $\{F_{p}(x)\}$ will become
$\{F_{\mu}(x)\}$ with $\mu=0,1,\cdots,p-1$
after a finite steps of the Collatz operation
as shown in Proposition 2. Hence, by using the
induction on the degree of the polynomials,
$\{F_{p}(x)\}$  for any $p$ are the Collatz
polynomials.

\vskip .3cm
In summary, the polynomials in representing
integers based on {a} binary numeral system
are introduced to explore
the Collatz conjecture, which seem
more convenient in the computation under the Collatz operation.
Especially, the polynomial structure
and its evolution under the Collatz
operation become more transparent, from which
the upper bound of the degree
of the polynomial after a finite steps of
the Collatz operation is estimated.
With this upper bound, it is shown
that the conjecture is
true in terms of the induction with
respect to the degree of the polynomials.

\section{acknowledgement}

{Support from the National Natural Science Foundation of China (11675071),
the  U. S. National Science Foundation (OIA-1738287 and ACI -1713690),
{U. S. Department of Energy (DE-SC0005248)}, the Southeastern Universities Research Association,
and the LSU--LNNU joint research
program (9961) is acknowledged.}



\end{document}